\documentclass{scrartcl}
\usepackage[utf8]{inputenc}
\usepackage[T1]{fontenc}

\usepackage{lmodern}
\usepackage{amsthm}
\usepackage{amsfonts}
\usepackage{mathtools}
\usepackage{stmaryrd}
\usepackage{MnSymbol}
\usepackage{xspace}
\usepackage{tikz-cd}
\usepackage[french=guillemets]{csquotes}
\MakeOuterQuote{"}
\usepackage[shortlabels]{enumitem}
\setlist{nosep}
\usepackage{xcolor}
\usepackage[style=alphabetic,backend=biber,sorting=nyt,giveninits=true]{biblatex}
\usepackage[colorlinks=true, breaklinks=true, urlcolor= black, linkcolor= black, citecolor= black, 
bookmarksopen=true,linktocpage=true,plainpages=false,pdfpagelabels]{hyperref}
\usepackage[capitalise]{cleveref}

\swapnumbers
\newtheorem{theorem}{Theorem}[section]
\newtheorem*{theorem*}{Theorem}
\newtheorem{claim}{Claim}[theorem]
\newtheorem{lemma}[theorem]{Lemma}
\newtheorem{corollary}[theorem]{Corollary}
\newtheorem{fact}[theorem]{Fact}
\newtheorem{proposition}[theorem]{Proposition}
\theoremstyle{definition}
\newtheorem{remark}[theorem]{Remark}

\newtheorem{definition}[theorem]{Definition}

\newtheorem{example}[theorem]{Example}

\newcounter{claimqed}
\newcommand\claimqed{\setcounter{claimqed}{0}%
\leavevmode\unskip\penalty9999
\hbox{}\nobreak\hfill\quad\hbox{$\square_{\text{Claim}}$}}
\newenvironment{claimproof}[1]{\par\noindent\textit{Proof.}\space#1\setcounter{claimqed}{1}}{\ifnum\value{claimqed}=1\claimqed\fi}

\newcommand\cC{\mathcal{C}}
\newcommand\cG{\mathcal{G}}
\newcommand\cX{X}
\newcommand\cY{Y}
\newcommand\A{A}
\newcommand\B{B}

\newcommand\D{\textup{(D)}\xspace}
\newcommand\I{\mathrm{I}}

\newcommand\Y{Y}
\newcommand\Z{Z}
\newcommand\K{\mathrm{K}}
\newcommand\V{V}
\newcommand\W{W}
\newcommand\f{f}

\newcommand\p{p}
\newcommand\q{q}

\newcommand\CC{\mathbb{C}}
\newcommand\RR{\mathbb{R}}
\newcommand\ZZ{\mathbb{Z}}
\newcommand\QQ{\mathbb{Q}}

\newcommand\m{\mathfrak{m}}

\newcommand\acl{\mathrm{acl}}
\newcommand\dcl{\mathrm{dcl}}
\newcommand\eq[1]{#1^\mathrm{eq}}
\newcommand\acleq{\eq{\acl}}
\newcommand\dcleq{\eq{\dcl}}

\renewcommand\k{\mathrm{k}}
\newcommand\res{\mathrm{res}}
\newcommand{\Stab}{\operatorname{Stab}}
\newcommand{\tp}{\operatorname{tp}}
\renewcommand\epsilon{\varepsilon}
\newcommand\Lin{\mathrm{Lin}} 
\newcommand\Int{\text{-}\mathrm{Int}}
\newcommand\kInt{\k\Int}
\newcommand\Aut{\mathrm{Aut}}
\newcommand\RV{\mathrm{RV}}
\newcommand\rv{\mathrm{rv}}
\newcommand\val{\mathrm{v}}
\newcommand\Gr{\mathrm{Gr}}
\newcommand\Cut{\mathrm{Cut}}
\newcommand\Mod{\mathrm{Mod}}
\renewcommand\mid{:}
\newcommand\pCut{\Cut^\star}
\newcommand\aCut{\Cut^{\star\star}}
\newcommand{\red}{\operatorname{red}}

\newcommand{\restr}[2]{\left.#1\right|_{#2}}

\newcommand\Hen{\mathrm{Hen}}
\newcommand\ACVF{\mathrm{ACVF}}
\newcommand\ACFA{\mathrm{ACFA}}

\newcommand\VFA{\mathrm{VFA}}

\newcommand{\MM}{\mathbb{M}}
\newcommand{\NTP}{\mathrm{NTP}}

\newcommand\RCVF{\mathrm{RCVF}}

\newcommand\ind{\downfree}
\newcommand\Linind[1][]{\ind^{\Lin_{#1}}}
\newcommand\alg[1]{#1^{\mathrm{a}}}
\newcommand\ur[1]{#1^{\mathrm{ur}}}
\newcommand\tensor{\otimes}
\newcommand\subsel{\preccurlyeq}
\newcommand\imp{\ensuremath{\Rightarrow}}

\renewcommand\L{\mathcal{L}}
\newcommand\G{\mathcal{G}}
\newcommand\Vg{\Gamma}
\renewcommand\O{\mathcal{O}}
\newcommand\code[1]{\ulcorner#1\urcorner}
\newcommand\qftpLin{\tp_{1}}

\newcommand\germ[2]{[#1]_{#2}}

\renewbibmacro{in:}{}
\bibliography{biblio}

\title{On residual domination and types orthogonal to the value group}

\author{Pablo Cubides Kovacsics\thanks{The author was supported by the research grant FAPA (INV-2023-158-3157) \emph{New methods in non-archimedean geometry and its applications}, funded by Universidad de los Andes.}
\and Silvain Rideau-Kikuchi
\and Mariana Vicar\'ia\thanks{The author was supported by the Humboldt Fundation (Postdoctoral Fellowship).}}

\begin{document}

\maketitle

\begin{abstract}
  We present a unifying framework of residual domination for (expansions of)
  henselian valued fields of equicharacteristic zero, encompassing some valued
  fields with operators. We show that the class of residually dominated types
  coincides with the types that are orthogonal to the value group, and with the
  class of types whose reduct to $\ACVF$ (the theory of algebraically closed
  valued fields with a non-trivial valuation) are generically stable. When the
  residue field is stable (resp. simple) we relate these equivalent notions to
  generic stability (resp. simplicity). Those results apply in particular to ultraproducts of $p$-adic fields and to the limit theory $\VFA_{0}$ of algebraically closed valued fields of characteristic $p$
  with the Frobenius automorphism (as $p$ tends to infinity).
  \renewcommand*{\thefootnote}{}
  \renewcommand*{\thefootnote}{\arabic{footnote}}
\end{abstract}

\section{Introduction}

One of the most celebrated results in recent model theory of valued fields is
Hrushovski-Loeser's characterization of the homotopy type of the Berkovich
analytification of quasi-projective varieties over complete (non-trivially
valued) non-archimedean valued fields \cite[Theorem $11.1.1$]{HL}. Inherent to
their approach, they introduced a model theoretic analogue of Berkovich's
analytification  of an algebraic variety $X$ as the  space of types concentrating
on $X$ which are orthogonal to the value group. One of the key features of that approach is the result proven in \cite{equivalence, HHM}:
in $\ACVF$ (the theory of algebraically closed non-trivially valued fields) the
following conditions are equivalent for a global $A$-invariant type $p$:
\begin{itemize}
\item $p$ is orthogonal to the value group,
\item $p$ is stably dominated,
\item $p$ is generically stable. 
\end{itemize}

The aim of the present paper is to extend such a theorem to arbitrary henselian
fields of equicharacteristic zero. Part of the difficulty lies in understanding
what should play the role of stability in such a generality. For example,  note
that, in the case of real closed valued fields, a generically stable type over a
model must be realized, and hence the previous equivalence is clearly false. It
turns out that, independently of what the stable part of the structure is or of
what generically stable types correspond to, the types which are orthogonal to
the value group coincide with those which are "controlled" by the residue field,
meaning that they are \emph{residually dominated} (see \cref{def resdom}). The
notion we introduce here encompasses notions of domination previously introduced
in \cite{has-ealy-mar, has-ealy-simon, Vic-ResDom} where they
are mostly studied over models.

Our contribution here is two-folded: first, we
extend such results over algebraically closed sets of imaginary elements and,
second, we show that they remain valid in various expansions of henselian valued
fields. In addition, we show that our notion of residual domination can be
evaluated by looking at the corresponding type in an algebraically closed
extension. Here is an informal version of the main contribution of this paper.

We refer the reader to \cref{equiv resdom} for the precise statement, to
\cref{propd} for the definition of property \D and to \cite{Vic-EIGp} for
background on bounded regular ordered abelian groups --- which are those with
countably many definable convex subgroups.

\begin{theorem*} Let $M$ be
an equicharacteristic zero henselian valued field such that the value group
$\Gamma(M)$ is \emph{either}:
\begin{itemize}
\item dense with property \D, or
\item a pure discrete ordered abelian group of bounded regular rank.
\end{itemize}
Then, the following are equivalent for a type $p$ of field points over an
algebraically closed set of imaginary elements:
\begin{enumerate}[(i)]
\item the type $p$ is residually dominated;
\item the quantifier free part of $p$, seen as a type in $\ACVF$, is residually dominated (equivalently orthogonal to \(\Gamma\), equivalently generically stable).
\end{enumerate}
\end{theorem*}

This theorem relies on descent results of \cite{descent} for stable domination
in arbitrary theories; but, as a corollary, we obtain base change for residual
domination (\cref{base change}) even when the residue field is not stable. A similar equivalence was proved independently by P.
Wang and D. Mutlu \cite{mutlu-wang} for invariant types over models.

We also extend the equivalence above to valued fields with operators (see
\cref{equiv resdom op}). This results applies in particular to valued difference
fields as in \cite{diffRV}, equicharacteristic zero \(\partial\)-henselian
fields with a monotone derivation \cite{Scanlon}, as well as models of
\(\Hen_{(0,0)}\) with generic derivations as in \cite{cubides-point} and
\cite{ForTer-GenDer}.



We then apply those results to show that residually dominated types inherit tame
behavior from the residue field. In particular, when the residue field is stable
we show that an type of field points is residually dominated if and only if it
is generically stable (see \cref{gen st equiv}). Similarly, under mild
hypotheses, if the residue field is simple, then a type of field points is residually dominated if and only if it is generically simple (see \cref{def:gen simple} and \cref{gen simp equiv}). Chernikov and Hils asked in \cite{VFANTP} whether there is a theory of simple domination analogous to the theory of stable domination as developed in \cite{HHM}. Our results around generically simple types can be seen as positive evidence of the existence of such a theory. In particular, these last results apply to ultraproducts of the $p$-adics and existentially closed multiplicative difference valued fields (see
\cref{pseudolocal and VFA}). One particularly striking example is the limit
theory $\VFA_{0}$ of algebraically closed valued fields of characteristic $p$
with the Frobenius automorphism (as $p$ tends to infinity).

\medskip

The main motivation for this work is the question of whether Hrushovski-Loeser's
ideas can be adapted to characterize the homotopy type of topological spaces
arising from types orthogonal to the value group. When the underlying theory is
the theory of real closed valued fields, one can view such spaces of types as
the Berkovich analytification of a semi-algebraic set as defined in
\cite{jell_etal}. To characterize the homotopy type of such spaces remains an
open question. It is also related to the wildly open question about the homotopy
type of spaces of real places. Our hope is to use the results of this paper
together with structural properties of such type spaces (for example strict
pro-definability, which was shown in \cite{beautiful}) to attack these questions.\\

\subsection*{Acknowledgments}
We want to thank M. Malliaris, Y. Dor and Y. Halevi for insightful
conversations.  The third author wishes to thank M. Malliaris also for
supporting her visit to the University of Chicago, where some of this research
took place; as well to A. Chernikov for the encouragement to pursue this
project and financially supporting the visit of the first author to UCLA.

\section{Preliminaries}

We refer the reader to \cite{TenZie} for an introduction to model theory. By
"definable" we mean "definable with parameters" unless we specify a language or
parameters in which case we mean "definable without (extra) parameters". We will
mostly abuse notations and write \(\L\) instead of \(\eq{\L}\).

Throughout this section, let \(M\) be a sufficiently saturated and homogeneous
structure.

\subsection{Stably embedded sets}

Let $A\subseteq \eq{M}$ and let $\cX$ and \(\cY\) be (ind-)$A$-definable sets.
We often identify (small) families of definable sets \((X_i)_{i\in I}\) with the
associated ind-definable set \(\bigcup_i X_i\).

\begin{definition}
\begin{enumerate}
\item The set \(\cX\) is \emph{stably embedded over $A$} if any \(M\)-definable
subset of \(\cX^n\), for any \(n\), is definable with parameters from $A \cup
\cX(M)$.
In that case, we denote by \(\eq{\cX}\) the family of all quotients of
\(A\)-definable subsets of \(\cX^n\), for any \(n\), by \(A\)-definable
equivalence relations.

\item The family \(\cX\) is \emph{$\cY$-internal} if there exists an
(ind-)\(M\)-definable subset \(Z\) of \(\cY\) and a surjective
(ind-)$M$-definable function $g : Z \rightarrow X$. We write \(\cY\Int_A\)
to denote the family of all $A$-definable $\cY$-internal sets.

\item The family \(\cX\) is \emph{almost $\cY$-internal} if there exists an
(ind-)\(M\)-definable subset \(Z\) of \(\cY\) and a surjective
(ind-)$M$-definable one-to-finite correspondence $g : Z \rightarrow X$.

\item The families \(\cX\) and \(\cY\) are \emph{orthogonal} if any definable
subset \(Z\) of \(\cX^n\times\cY^m\), for any \(n\) and \(m\), is a finite union
of squares; \emph{i.e.} there are finitely many \(M\)-definable sets $V_{\ell}$
of \(\cX\) and \(W_\ell\) of \(\cY\) such that \(Z=\bigcup_{\ell} V_{\ell}
\times W_{\ell}\).

\item Assume \(\cX\) is stably embedded over \(A\). A global \(A\)-invariant type
\(p\) is \emph{orthogonal to \(\cX\)} if, for any \(a\models p\), we have
$\eq{\cX}(Ma)\subseteq\eq{M}$.
\end{enumerate}
\end{definition}

When \(\cX = \varinjlim_i X_i\) is stably embedded over \(A\), we will often
consider it as a structure whose sorts are the \(X_i\) with the full
\(A\)-induced structure on products of sorts.\\

We now assume that \(\cX\) and \(\cY\) are stably embedded over \(A\). The
following statements are folklore (\emph{e.g.} by \cite[Lemma~2.5.(c)]{CDHJR} for
related statements).

\begin{lemma}
\label{st emb prop}
Let $A \subseteq B\subseteq\eq{M}$ and let $a$ be a tuple in $\cX$.
\begin{enumerate}
\item If $\cX$ codes finite sets then
$\cX(\acleq(B))\subseteq \acleq(A \cX(B))$.

\item We have $\tp(a/A \eq{\cX}(B)) \vdash \tp(a/B)$.
\end{enumerate}
\end{lemma}

In particular, $\eq{\cX}(\acleq(B))\subseteq \acleq(A \eq{\cX}(B))$ and, if
$\cX$ eliminates imaginaries, we have $\tp(a/A\cX(B)) \vdash \tp(a/B)$.

\begin{proof}
Regarding the first statement, let $d \in \cX(\acleq(B))$ and let $F$ be
the finite set of conjugates of $d$ over $B$. Since $\cX$ codes finite
sets, there is $e \in \cX$ interdefinable over $A$ with $\code{F}$. We
have $\ulcorner F \urcorner \in \dcleq(B)$ so $e \in \cX(B)$. Since
\(d\in F\), we have $d \in \acleq(A \ulcorner{F}\urcorner)\subseteq
\acleq(A e)\subseteq \acleq(A \cX(B))$.

As for the second statement, let \(Y\in \tp(a/B)\). As \(\cX\) is stably
embedded \(Y\) is \(\cX(M)\)-definable so \(\code{Y}\in \eq{\cX}(B)\).
Consequently, $\tp(a/A \eq{\cX}(B)) \vdash \tp(a/B)$.
\end{proof}

\begin{fact}
\label{type game}
Let \(A\subseteq M\) and let \(a,b,c\) be tuples in \(M\). The following are
equivalent:
\begin{enumerate}
\item \(\tp(b/Ac) \vdash \tp(b/Aa)\);
\item \(\tp(ac/A)\cup\tp(bc/A)\vdash \tp(abc/A)\);
\item \(\tp(ac/A)\cup\tp(c/Ab)\vdash \tp(ac/Ab)\).
\end{enumerate}
\end{fact}

Let us now consider properties of orthogonality.

\begin{lemma}
The sets \(\cX\) and \(\cY\) are orthogonal if and only if for every
\(M\)-definable function \(Z\subseteq X^n\), every $\f: Z \rightarrow \eq{\cY}$ has finite
image.
\end{lemma}

\begin{proof}
First assume that the families are orthogonal and let $\f: Z \rightarrow
\eq{\cY}$ be definable. For every \(x\in \cX\), let \(F_x\subseteq \cY^m\)
denote elements of the equivalence class \(f(x)\) and let $F=\{(x,y) \in Z
\times \cY^m \mid y \in F_x\}$. By orthogonality, $F=\bigcup_{\ell\leq k} \V_{\ell}
\times \W_{\ell}$ for some definable $\V_{\ell} \subseteq \cX^n$ and $\W_{\ell}
\subseteq \cY^m$. We may assume that the \(\V_\ell\) form a partition of
\(\cX\). Then \(\f(\cX) = \{\ulcorner \W_{\ell}\urcorner \mid \ell\leq k\}\).

Conversely, let $F \subseteq \cX^n \times \cY^m$ be definable. For every \(x\in
\cX\), let \(F_x =\{ y \in \cY^m \mid (x,y) \in F\}\). Let $Z$ be the projection of $F$ unto $X^n$ and let $\f: \Z
\rightarrow \eq{\cY}$ send $x$ to $\ulcorner F_{x} \urcorner$. Since $\f$ has
finite image, there are finitely many definable sets $ \W_{1},\dots, \W_{k}
\subseteq \Y$ such that for all $x \in \V$, we have $F_{x}=\W_{\ell}$ for some
$\ell$. Let $\V_{\ell}=\{x \in \cX^n \mid F_{x}=\W_{\ell}\}$, then $F=
\bigcup_{\ell\leq k} \V_{\ell} \times \W_{\ell}$.
\end{proof}

In particular, if \(\cX\) and \(\cY\) are orthogonal, then any invariant type
concentrating on \(\cX\) is orthogonal to \(\cY\).

\begin{remark}
An invariant type \(p\) is orthogonal to \(\cX\) if and only if for every
definable function \(f: \V \to \eq{\cX}\), with \(p\) concentrating on \(V\),
the type \(f(p)\) is realized. Indeed, if \(p\) is orthogonal to \(\cX\), then
for any \(a\models p\), we have \(f(a) \in \eq{\cX}(Ma) \subseteq\eq{M}\) and
hence \(f(p) = \tp(f(a)/M)\) is realized. Conversely, if every such type is
realized, then \(f(a) \in \eq{M}\).
\end{remark}

Let us conclude this section with considerations on forking. 

\begin{lemma}
\label{fork st emb}
Let \(a\) be a tuple of elements from \(\cX\), let $B\subseteq C
\subseteq \eq{M}$ contain $A$. Then 
\begin{equation*}
a \ind_{B} C \text{ if and only if } a \ind_{A \eq{\cX}(B)} \eq{\cX}(C). 
\end{equation*}
\end{lemma}

In particular, if $\cX$ eliminates imaginaries, we have
\begin{equation*}
c \ind_{B} d \text{ if and only if } c \ind_{A \cX(B)} \cX(B d). 
\end{equation*}

\begin{proof}
We have $\tp(a/A \eq{\cX}(C)) \vdash \tp(a/ C)$ by \cref{st emb prop}.2.
If there is a definable set \(Y\in \tp(a/ C)\) that forks over \(B\), then there
is definable set in \(\tp(a/ A\eq{\cX}(C))\) contained in \(Y\) which
therefore forks over \(A \eq{\cX}(B)\).

Conversely, if $Y\in \tp(a/ A  \eq{\cX}(C))$ forks over
$A\eq{\cX}(B)$, then \(Y = \bigcup_i Z_i\) where the \(Z_i\) are
\(M\)-definable and divide over \(A\eq{\cX}(B)\). As \(\cX\) is stably
embedded, any \(A\eq{\cX}(B)\)-indiscernible sequence containing
\(\code{Z_i}\in \eq{\cX}(M)\) is also \(B\)-indiscernible. So \(Z_i\) divides
over \(B\) and \(Y\) forks over \(B\).
\end{proof}

\begin{remark}[{\cite[Proposition~2.12]{3dividing}}]
\label{forking acl}
Forking independence preserves algebraic closure in any theory; that is for any
$A,B,C \subseteq M$, we have
\begin{align*}
A \ind_{B} C&\Leftrightarrow \acl(AB) \ind_{B} C\\
&\Leftrightarrow A \ind_{B} \acl(BC)\\
&\Leftrightarrow A \ind_{\acl(B)} C.
\end{align*}
\end{remark}

\subsection{Domination}

Let \(M\) be some sufficiently saturated and homogeneous structure, let
\(A\subseteq\eq{M}\), let \(\cX\) be a stably embedded ind-\(A\)-definable set.

\begin{definition}
Let $p$ be a type over $A$ and \(f\) be an pro-\(A\)-definable map into \(\eq{\cX}\) which is defined at $p$. We say that:
\begin{itemize}
\item a type $p$ over $A$ is dominated by \(\cX\) via \(f\) if for any realization $a\models p$ and
every \(B\supseteq A\) such that \(f(a) \ind_A B\), we have \(\tp(B/Af(a))
\vdash \tp(B/Aa)\);

\item a type $p$ over $A$ is dominated by \(\cX\) if it is so via any pro-\(A\)-definable function $f$ sending $x$ to \(\eq{\cX}(\dcleq(Ax))\);
\item a type $p$ over $M$ is dominated by \(\cX\) over \(A\) if it is \(A\)-invariant
and $\restr{p}{A}$ is dominated by \(\cX\).
\end{itemize}
\end{definition}

\begin{remark}
\label{models suffice dom} 
To check that a type \(p\) over \(A\) is dominated by \(X\) via \(f\), it suffices that for
some \(a\models p\) and some sufficiently saturated and homogeneous \(N\subsel M\) containing \(A\) with
\(f(a) \ind_A N\), we have \(\tp(N/Af(a)) \vdash \tp(N/Aa)\). Indeed, let \(a'\models p\) and \(B'\supseteq A\). To determine \(\tp(B'/Aa')\) it suffices to consider finite tuples \(b'\) in \(B'\). Let \(b\) be such that \(ab \equiv_A a'b'\). Up to \(A\)-isomorphism, we may assume that \(b\subseteq N\). By \cref{type game}, we have \(\tp(Nf(a)/A) \vdash \tp(N/Aa)\) and thus \(\tp(b/Af(a))\vdash
\tp(b/Aa)\). It follows that \(\tp(b'/Af(a'))\vdash
\tp(b'/Aa')\). 
\end{remark}

\begin{remark}
\label{dom fork}
Assume that \(\cX\) is stable and let \(a\) be a tuple in \(M\).
\begin{enumerate}
\item If \(\tp(a/A)\) is dominated by \(\cX\) via some \(f\), then the non
forking (equivalently \(\acleq(A)\)-definable) extensions of \(\tp(a/A)\) over
\(A\) correspond to the non forking (equivalently \(\acleq(A)\)-definable)
extensions of \(\tp(f(a)/A)\) --- see \cite[Lemma~3.12]{HHM}. In particular, for
every \(B\supseteq A\), if \(f(a)\ind_A B\), then \(a\ind_A B\).
It follows that \(\tp(a/A)\) is stationary if and only if \(\tp(f(a)/A)\),
equivalently \(\tp(a/A)\), has an \(A\)-invariant extension.

\item Conversely, if \(\tp(a/A)\) is stationary, then it is dominated by \(\cX\)
via \(f\) if and only if \(f(a)\ind_A B\) implies \(a\ind_A B\), for every
\(B\supseteq A\).
\end{enumerate}
\end{remark}




\begin{definition}
\label{def:gen st}
A type $\p(x)$ over \(M\) is \emph{generically stable over $\A$} if it is
$\A$-invariant and for any ordinal $\alpha \geq \omega$, for any Morley sequence
$(a_{i})_{i < \alpha}$ of $\p$ over $\A$ --- i.e $a_{i} \vDash \restr{\p}{\A
a_{<i}}$ --- and for any $\L(M)$ formula $\phi(x)$, the set $\{ i < \alpha \mid
M \vDash \phi(a_{i})\}$ is finite or co-finite.
\end{definition}

\begin{remark}
\label{resdom gen st}
Assume \(\cX\) is stable and let \(p\) be a stationary type over \(A\) which is dominated by \(X\). Then its non forking extension \(q\) is generically stable and so is \(q^{\tensor\omega}\) --- see \cite[Fact~2.45]{descent}.
\end{remark}

\begin{theorem}
\label{descent intdom}
Let $p$ be global $A$-invariant type. Assume that it is dominated by some
$A$-definable stable set $\cX$. Then \(\restr{p}{A}\) is dominated by
$\cX\Int_{A}$.
\end{theorem}

\begin{proof}
This is the content of \cite[Theorem~4.3]{descent}. The hypothesis are granted by \cref{resdom gen st}. 
\end{proof}

\begin{definition}
Let \(p\) be an \(A\)-definable type over \(M\). Two definable functions \(f\)
and \(g\) defined at realizations of \(p\) have the same \emph{germ at \(p\)},
if \(p(x)\vdash f(x) = g(x)\). We write \(\germ{f}{p} \in\eq{M}\) for the class
of functions with the same germ at \(p\) as \(f\) varies in a
\(\emptyset\)-definable family. The germ of \(f\) at \(p\) is \emph{strong}
(over \(A\)) if, for any \(a\models \restr{p}{A}\), we have
\(f(a)\in \dcleq(A,\germ{f}{p},a)\).
\end{definition}

\begin{remark}
\label{strong germs}
If \(\cX\) is stable and \(p\) is stably dominated over \(A\)
(\cite[Theorem~6.3]{HHM}), or more generally if \(p\) is generically stable
(\cite[Theorem~2.2]{AdlCasPil}), then all germs at \(p\) are strong.
\end{remark}

\subsection{Equicharacteristic zero henselian fields}

We refer the reader to \cite{EngPre} for background on valued
fields. If \((F,\val)\) is a valued field, we write \(\Vg(F)\) for the
\emph{value group} \(\val(F^\times)\), \(\O(F)\) for its valuation ring,
\(\m(F)\) for its maximal ideal and \(\k(F)\) for its residue ring and \(\res :
\O(F)\to\k(F)\) the residue map. We denote by \(\RV_F^\times(F)\) the group
\(F^\times/(1+\m(F))\) of \emph{leading terms} and \(\rv : F\to\RV(F) =
\RV^\times(F) \cup\{0\}\) the canonical projection, extended by \(\rv(0) = 0\).
Observe that it sits in the short exact sequence
\[1\to\k^\times(F)\to\RV^\times(F) \to \Gamma(F) \to 0.\]

We consider valued fields in the language \(\L_\K\) with one sort \(\K\) for the
valued field, with the ring structure and a predicate for
\(\val(x)\leq\val(y)\). The residue field, the value group and the leading terms
are then interpretable and when we need to consider such sets we will be working
in the language \(\eq{\L_\K}\) where we add all \(\L_\K\)-interpretable
sets.

If $S$ is any of the sort of \(\eq{\L_\K}\), by an \emph{$S$-expansion} we mean
a language extending $\eq{\L_\K}$ which adds functions, relations and constant
symbols to $S$.

A central result in model theory of henselian valued fields in equicharacteristic
zero is relative quantifier elimination down to \(\RV\). We refer the reader to
\cite[Appendix~A]{Rid-ADF} for background on relative quantifier elimination
(and expansions).

\begin{theorem}[{\cite[Theorem~2.1]{Del-These},\cite[Theorem~B]{Bas}}]
\label{EQ RV}
Let \(M\) be an \(\RV\)-expansion of an equicharacteristic zero henselian valued
field and let $A\leq \K(M)$ is a sub-ring. Every $A$-definable subset of \(\K^n\)
is of the form
\[\{x \in \K^n\mid \rv(P(x)) \in X\}\] where \(P\) is an \(m\)-tuple in \(A[x]\)
and \(X\subseteq \RV^{m}\) is \(\rv(A)\)-definable.
\end{theorem}

When the residue field is algebraically closed quantifier elimination can be
further simplified:

\begin{fact}[see {\cite[Corollary~2.33]{Vic-EI}}]
\label{EQ unramified}
Let \(M\) be a \(\Vg\)-expansion of an equicharacteristic zero henselian valued
field with algebraically closed residue field and let  $A\subseteq \K(M)$ be a
sub-ring. Every \(A\)-definable subset of \(\K^n\) is of the form
\[\{x\in \K^n\times\Vg^m\mid \val(P(x)) \in X\}\] where \(P\) is an \(m\)-tuple
in \(A[x]\) and \(X\subseteq \Gamma^{m}\) is \(\val(A)\)-definable in \(\Vg\).
\end{fact}

As a byproduct of quantifier elimination, we have the following:

\begin{fact}[{\cite[Corollary~2.25 and 2.24]{vdD-notes}}]
\label{orth Vg k}
Let \(M\) be an equicharacteristic zero henselian valued field. Then
\begin{itemize}
\item \(\k\) is a stably embedded pure field;
\item \(\Vg\) is a stably embedded pure ordered abelian group;
\item \(\k\) and \(\Vg\) are orthogonal.
\end{itemize}
\end{fact}

In \cite{RKVic-EIAKE}, to classify interpretable sets, the authors introduce the
(generalized) geometric sorts --- building on the geometric sorts of Haskell,
Hrushovski and Macpherson \cite{HHM-EI}. Let us briefly recall their definition.
For an in depth discussion of the geometric sorts, we refer the reader to
\cite[Section~3]{RKVic-EIAKE}.

Let \(M\) be a valued field. We fix an (ind-)\(\emptyset\)-definable family
\(\cC = (\cC_c)_{c\in\Cut}\) of cuts in \(\Gamma\) such that any \(M\)-definable
cut is of the of form \(\cC_c\) for some unique \(c \in \Cut(M)\). We let
$\pCut$ be the subset of proper cuts. For every \(c\in\Cut\), let \(\I_c\)
denote the \(\mathcal{O}\)-submodule \(\{x\in \K\mid\val(x) \in \cC_c\}\). Let
\(\B_n(\K)\) denote the group of $n \times n$ upper triangular invertible
matrices with coefficients in \(\K\). It acts naturally on (definable)
\(\O\)-submodules of \(\K^n\).

\begin{definition}
\label{geom sorts}
\begin{enumerate}
\item For every \(n\)-tuple \(c\) in \(\Cut^\star\), let \(\Lambda_c\) be the
\(\O\)-module \(\sum_{i\leq n} \I_{c_i}e_i\), where \((e_i)_i\) is the canonical
basis of \(\K^n\).
\item Let \(\Mod_c := \B_n/\Stab(\Lambda_c)\). We identify any element \(A
\Stab(\Lambda_c)\) of \(\Mod_c\) with the definable \(\O\)-module \(A \Lambda_c
\subseteq \K^n\).
\item Let \(\aCut\) denote the cuts that are distinct from \(+\infty\) and
\(-\infty\) as well as, if \(\Gamma(M)\) is dense any \(\gamma^+\) with
\(\gamma\in\Vg\).
\item For every definable \(\mathcal{O}\)-module \(R\), let \(\red(R)\) denote
the \(\k\)-vector space \(R/\m R\).
\item Let \(\Gr := \bigsqcup_{c\subseteq \aCut} \Mod_c\) and let \(\Lin :=
\bigsqcup_{R\in \Gr} \red(R)\). The geometric sorts \(\G\) consists of
\(\K\cup\Gr\cup \Lin\).
\end{enumerate}
\end{definition}

The main raison d'être of the geometric sorts is the following:

\begin{proposition}[{\cite[Corollary~3.21]{RKVic-EIAKE}}]
\label{code module}
Let \(M\) be an \(\RV\)-expansion of an equicharacteristic zero henselian valued
field. Then any definable \(\O\)-submodule \(R\subseteq \K^n\) is coded in
\(\G\).
\end{proposition}

\subsection{The structure \texorpdfstring{\(\Lin_A\)}{LinA}}

Let \(M\) be an expansion of a valued field and let $\A \subseteq \eq{M}$. We
assume \(\k\) to be stably embedded.

\begin{definition}
\label{LinA}
We consider the ind-\(A\)-definable set
\[\Lin_\A := \bigsqcup_{R\in\Gr(\acleq(\A))}\red(R).\]
\end{definition}

For every \(R\in \Gr\), \(\red(R)\) is a finite dimensional $\k$-vector space.

\begin{remark}
\label{rem: linfacts}
The set ind-definable set \(\Lin_A\) is stably embedded as it suffices to name
bases to identify each $\red(R)$ with \(\k^m\) for some $m$ (see \cite[Lemma
2.5.18]{HilRK-EIAKE}). We consider \(\Lin_A\) as a structure with a sort for
\(\k\) with its \(A\)-induced structure and a sort for each \(\red(R)\) with its
vector space structure. It is a \(\k\)-linear structure in the sense of
\cite{Hru-Groupoid}. It is closed under tensors, duals and has flags by
\cite[Proposition~3.19]{RKVic-EIAKE}. Moreover, If $\k$ is stable (resp. simple)
then $\Lin_{A}$ is also stable (resp. simple), as each of the finite dimensional
$\k$-vector spaces \(\red(R)\) is internal to $\k$.\\

If \(\Gr(\acleq(A))\subseteq A\), the \(A\)-induced structure on \(\Lin_A\) is a
definable expansion of this linear structure.
\end{remark}

\begin{proposition}[{\cite[Lemma 5.6]{Hru-Groupoid}}]
\label{LinA EI}
If \(\k\) is a pure algebraically closed field and if \(\Gr(\acleq(A))\subseteq
A\), then \(\Lin_A\) eliminates imaginaries.
\end{proposition}

Our interest in \(\Lin_A\) is that, under mild hypotheses, it encompasses all
\(A\)-definable \(\k\)-internal sets.

\begin{definition}
\label{propd}
We say that an (expansion of an) ordered set \(G\) satisfies Property $\D$ if
for every finite set of formulas $\Delta(x,y)$ containing the formula $x<
y_{0}$, any \(A = \acleq(A)\subseteq \eq{G}\) and any \(\Delta\)-type $p(x)$
which is $\A$-definable, there is an $\A$-definable complete type $\q(x)$
containing \(\p\). 
\end{definition}

\begin{remark}
\label{propertyD}
Property $\D$ is a stronger than density of definable types. The third author
proved that it holds in pure ordered abelian groups of bounded regular rank (see
the second half of the proof of \cite[Theorem~5.3]{Vic-EI}). Note that it holds
trivially for any $o$-minimal expansion of an ordered abelian group. 
\end{remark}

\begin{proposition}
\label{char kint} 
Let \(M\) be a \(\Gamma\)-expansion of an equicharacteristic zero henselian
field whose residue field is algebraically closed and whose value group is dense
and satisfies property \D. Let \(A\subseteq \eq{M}\) and \(X\) be
\(A\)-interpretable. The following statements are equivalent:
\begin{enumerate}
\item \(X\) is \(\k\)-internal;
\item \(X\) is almost \(\k\)-internal;
\item \(X\) is orthogonal to \(\Gamma\);
\item there is an \(A\)-definable finite-to-one map from \(X\) into \(\Lin_A\).
\end{enumerate}
\end{proposition}

\begin{proof}
By \cite[Corollary~6.7]{RKVic-EIAKE}, the first three statements are equivalent to 
\begin{itemize}
\item[4'.] \(X \subseteq \dcleq(\acleq(A),\Lin_A(M))\).
\end{itemize}
Note that statement 4 implies statement 2 because $\Lin_A$ is $\k$-internal. Hence, it suffices to show that statement 4' implies 4. If statement 4' holds, then there is an \(\acleq(A)\)-definable function \(f\)
from \(\Lin_A\) to \(X\). Since \(\k\) is algebraically closed, \(\Lin_A\)
eliminates imaginaries (see \cite[Remark~3.20.2]{RKVic-EIAKE}) and hence for
every \(e\in X\), the set \(f^{-1}(e)\) is coded by some tuple \(\eta\) in
\(\Lin_A(\acleq(A)e)\). Let \(\overline{\eta}\) be the code in \(\Lin_A\) of the
finite orbit of \(\eta\) over \(Ae\). Then \(\overline{\eta}\in \dcleq(Ae)\) and
\(e\in \acleq(A\overline{\eta})\).
\end{proof}

\section{Setup}\label{setup}

Let us fix some notation for the rest of the paper. Let \(M\) be a sufficiently
saturated and homogeneous expansion of an equicharacteristic zero henselian
valued field in some language \(\L\). We assume that \(M\) admits a spherically
complete elementary extension and that \(\k\) and \(\Vg\) are stably embedded
and orthogonal. This is the case, for example, if there is no expansion
(\cref{orth Vg k}). We further assume that \(\Vg(M)\) is \emph{either}:
\begin{itemize}
\item dense with property \D, or
\item a pure discrete ordered abelian group with bounded regular rank --- in
which case we include a constant \(\pi\) for a uniformizer in \(\L\).
\end{itemize}
Let \(\MM\) be a sufficiently saturated and homogeneous elementary extension of
\(M\). We will also need to consider various (non elementary) extensions of
\(M\).
\begin{itemize}
\item Let \(M_0 = \alg{M}\) be the algebraic closure of \(M\) considered as a
structure in the language \(\L_0 =\L_\K\). Note that \(M_0\) eliminates
quantifiers (\emph{e.g.} by \cref{EQ unramified} and quantifier elimination in
divisible ordered abelian groups). Let \(\MM_0\) be a sufficiently saturated and
homogeneous extension of \(M_0\). By quantifier elimination, we may assume that
\(\alg{\MM} \subsel \MM_0\).

\item Let \(M_1 = \ur{M}\) be the maximal algebraic unramified extension of
\(M\). We consider it as a \(\Vg\)-expansion of \(\L_\K\) by all
\(\L\)-definable subsets of \(\Vg\) in \(M\) --- recall that \(\Vg(M_1) =
\Vg(M_2)\). By \cref{EQ unramified}, the structure \(M_1\) eliminates
quantifiers. Let \(\MM_1\subseteq \MM_0\) be  a sufficiently saturated and
homogeneous extension of \(M_1\). By quantifier elimination, we may, once again,
consider that \(\ur{\MM} \subsel \MM_1\).

\item If \(\Vg(M)\) is discrete, let \(M_{1'} = M_1[1/\pi^{1/\infty}]\) the
extension obtained by adding \(n\)-th roots of \(\pi\) for all \(n>0\). We
consider \(M_{1'}\) as a structure in the \(\Vg\)-expansion \(\L_{1'}\) of
\(\L_\K\) by predicates for all definable convex subgroups \(\Delta\leq
\Vg(M_{1'})\) and the language of Presburger arithmetic for all quotients
\(\Vg/\Delta\). As \(\Vg(M_{1'})\) has bounded regular rank
(\cite[Lemma~5.10]{RKVic-EIAKE}), it eliminates quantifiers (see
\cite{Vic-EIGp}) and hence so does \(M_{1'}\).

Let \(\MM_{1'}\subseteq\MM_0\) be a sufficiently saturated and homogeneous
elementary extension of \(M_{1'}\). Note that the structure induced by
\(M_{1'}\) and \(\MM_1[1/\pi^{1/\infty}]\) on \(\Vg(M_1)\) coincide; see
\cite[Lemma~5.10]{RKVic-EIAKE} for a description of all definable convex
subgroups of \(M_1\). So we may assume that \(\MM_1[1/\pi^{1/\infty}]\subsel
\MM_{1'}\).
\end{itemize}

If \(\Vg(M)\) is dense, by convention we take \(M_{1'} = M_1\) and \(\MM_{1'} =
\MM_1\). For \(\epsilon\in\{0,1,1'\}\), we let indices \(_\epsilon\)
indicate that we consider the structure \(\MM_\epsilon\). For
example, \(\dcleq_\epsilon\) denotes the definable closure in
\(\eq{\MM}_\epsilon\).
Similarly, given a sort \(S\) of \(\eq{\L}_\epsilon\) and \(A\subseteq
\eq{\MM}_\epsilon\), we write \(S_\epsilon(A)\) for \(\dcleq_\epsilon(A)\cap
S(\eq{M}_\epsilon)\).

\begin{example}
\begin{enumerate}
    \item Let $M$ be a model of $\RCVF$. We have $M_1\simeq M_0$.  
    \item Let $M=\RR(\!(t)\!)$. We have $M_1\simeq \CC(\!(t)\!)\cap M_0$ and we
    have $M_{1'}\simeq M_0$.
    \item Let $(\ZZ^2, +, 0, \leqslant_\mathrm{lex})$, let
    $\Gamma'=(\ZZ\times\QQ, +, 0, \leqslant_\mathrm{lex})$ and let \(M =
    \RR(\!(t^\Gamma)\!)\). Then we have $M_1\simeq\CC(\!(t^\Gamma)\!)\cap M_0$.
    Taking $\pi=t$, we have $M_{1'}\simeq \CC(\!(t^{\Gamma'})\!)\cap M_0$.
    \end{enumerate}
\end{example}

We conclude this section with a comparison of the types in the various
structure, when they are orthogonal to \(\Vg\):

\begin{lemma}
\label{orth Vg ext ACVF}
Let \(A\subseteq \eq{M}\) contain \(\dcleq(A)\cap\G_0\) and let \(a\) be a tuple
in \(\K(\MM)\) such that \(\tp_0(a/M)\) is \(\L(A)\)-invariant and
\(\val(M(a))=\val(M)\). Let \(\epsilon \in \{0,1,1'\}\). Then 
\begin{enumerate}
\item \(\tp_0(a/M) \vdash \tp_\epsilon(a/M_\epsilon)\);
\item \(\tp_\epsilon(a/M_\epsilon)\) is \(\L_\epsilon(\cG_0(A))\)-definable
and orthogonal to \(\Gamma\).
\end{enumerate}
\end{lemma}

\begin{proof}
We roughly follow the proof of \cite[Proposition~4.3]{RKVic-EIAKE}.

By \cite[Lemma 3.3.7]{HilRK-EIAKE}, the type \(\tp_0(a/M)\) has a unique extension to
\(M_0\), so \(\tp_0(a/M)\vdash \tp_0(a/M_0)\). Moreover, since \(\val(M(a)) =
\val(M)\), we have
\begin{equation}\label{Vg growth}
\val(M_\epsilon(a))\subseteq (\QQ\tensor \val(M(a)))\cap \val(\MM_\epsilon)
= \val(M_\epsilon).
\end{equation}
It follows from quantifier elimination in \(M_1\) (\cref{EQ unramified}) that
\(\tp_0(a/M_1) \vdash \tp_1(a/M_1)\); and hence that \(\tp_0(a/M)\vdash
\tp_1(a/M_1)\). Similarly, if \(\Gamma(M)\) is discrete, we have \(\tp_0(a/M)
\vdash \tp_{1'}(a/M_{1'})\).

Let us now consider definability of \(\tp_0(a/M_0)\). Growing \(M\), we may
assume that it is maximally complete. For every integer \(d\geq 0\), let $V_{d}
\simeq \K^\ell$ be the space of polynomials in $\K(M)[x]_{\leq d}$ of degree
less or equal than $d$ (in each \(x_i\)). It comes with a valuation \(v\)
defined by \(v(P)\leq v(Q)\) if \(\val(P(a))\leq\val(Q(a))\). Since \(M\) is
maximally complete, this vector space admits a separated basis \((P_i)_i\) (see
\cite[Proposition~3.3]{RKVic-EIAKE}). As \(\val(M(a)) = \val(M)\), we may
assume \(v(P_i) = 0\).

By \cite[Claim~3.3.5]{HilRK-EIAKE}, the family \((P_i)_{i}\) is also a separated
basis of \(V_d^0 = \K(M_0)[x]_{\leq d}\) for the valuation $v_{0}$ defined by
\(v_0(P)\leq v_{0}(Q)\) if \(\val(P(a))\leq\val(Q(a))\). So \(v_0\) is coded by
the \(\O\)-lattice \(\{ P \in V_{d}^0\mid v_0(P) \geq 0\}\) which we identify
with its intersection with \(M\) (\emph{cf.} \cite[Remark 4.2]{RKVic-EIAKE}). By
\cref{code module}, this lattice is \(\L_0(\G_0(A))\)-definable and hence so is
\(\tp_0(a/M_0)\).

Definability of \(\tp_\epsilon(a/M_\epsilon)\) follows and orthogonality to
\(\Gamma\) is a consequence of \cref{Vg growth}.
\end{proof}

\section{Residual domination}
Let \(A\subseteq\eq{M}\) be a small set of parameters. If \(B\subseteq\eq{M}\), we write \(\Lin_A(B)\) for
\(\Lin_A\cap\dcleq(AB)\).

The structure \(\alg{\k}\tensor\Lin_A\) is the linear
structure over the (pure) algebraically closed field \(\alg{\k}\) obtained by
tensoring each sort of \(\Lin_A\) with \(\alg{\k}\). It is a stable structure
that eliminates imaginaries (\cref{LinA EI}). We consider it as an
ind-\(\G(A)\)-definable structure in \(M_1\).

\begin{remark}
\label{dclka LinA}
Note that \(\k(M)\) is perfect and it contains bases for every \(\red(R)\) with
\(R\in\Gr(\acleq(A))\). So we have \(\dcleq_{1}(M)\cap\alg{\k}\tensor\Lin_A
\subseteq \Lin_A(M)\).

By stable embeddedness and elimination of imaginaries in (\cref{LinA EI}), it follows that if \(B\subseteq\eq{M}\) containing \(\Lin_A(\acleq(A))\), any \(\L(B)\)-definable subset of (a power of) \(\Lin_A\) is \(\Lin_A(B)\)-definable.
\end{remark}

\begin{definition}
Let \(a\) be a tuple in \(\Lin_A(M)\) and let \(B\subseteq C\subseteq\eq{M}\).
We write \(a\Linind_B C\) if \(a\) does not fork with \(\Lin_A(\acleq(A)C)\) over
\(\Lin_A(\acleq(A)B)\) in the structure \(\alg{\k}\tensor\Lin_A\).
\end{definition} 

Note that, by \cref{char Linind}.2, the relation \(\Linind\) does not depend on the choice of \(A\).

\begin{lemma}
\label{char Linind}
Let \(B \subseteq \eq{M}\) contain \(\Lin_A(\acleq(AB))\). Let \(a\) be a tuple
in \(\Lin_A(M)\).
\begin{enumerate}
\item The type \(\qftpLin(a/B)\) is \(\L_1(\Lin_A(B))\)-definable, stationary and its non forking extension is consistent with \(\tp(a/B)\).
\item Let \(C\subseteq \eq{M}\) contain \(B\). The following are equivalent:
\begin{enumerate}[(i)]
\item \(a\Linind_{B} C\);
\item \(\qftpLin(a/C)\) admits an \(\L(B)\)-invariant extension
\(\qftpLin(c/M)\) which is finitely satisfiable in \(M\);
\item \(a\) does not fork with \(C\) over \(B\) in \(M_1\).
\end{enumerate}
\end{enumerate}
\end{lemma}





\begin{proof}
Let us prove the first assertion. By \cite[Lemma~2.7]{Pil-GeomSta}, there exists
\(c\equiv_B a\) with \(\qftpLin(c/M)\) definable over \(\acl(AB)\). By
\cite[Lemma~3.3.7]{HilRK-EIAKE}, we have \(\qftpLin(c/M) \vdash \qftpLin(c/M_1)
:= p\). By \cref{dclka LinA},
the type \(p\) is coded, over \(\Gr(\acleq(A))\), by some tuple \(e\) in
\(\dcleq_{1}(\Lin_A(M)) \cap\alg{\k}\tensor\Lin_A\subseteq \Lin_A(M)\). As \(e\)
is fixed by every \(\sigma\in\Aut(M/\acleq(AB))\), we have \(e\in
\Lin_A(\acleq(AB))\subseteq \Lin_A(B)\) and thus \(p\) is \(\L_1(\Lin_A(B))\)-definable and stationary over \(B\).

Let us now consider the second assertion. We have \(a\Linind_{B} C\) if and only
if \(a\models \restr{p}{\Lin_A(C)}\) --- equivalently \(a\models \restr{p}{C}\), by \cref{dclka LinA} --- if and only if \(a\) does not fork with C over \(B\) in \(M_1\). Therefore it suffices to show that any
\(\L(B)\)-invariant type \(\qftpLin(c/M)\) consistent with \(\tp(a/B)\) is equal
to \(\restr{p}{M}\). Let \(c\equiv_B a\) such that \(c\models q\). As above, the
unique extension of \(q\) to \(M_1\) is stationary over \(B\) and therefore
coincides with \(p\).
\end{proof}

\begin{definition}
\label{def resdom}
Let $p$ be a type over $A$ (in \(\eq{M}\)) and let \(f\) be a pro-\(A\)-definable map
into a power of \(\Lin_A\) which is defined at $p$. We say that:
\begin{itemize}
\item $p$ is residually dominated via \(f\) if for every realization $a\models p$ and for every
\(B\subseteq\eq{M}\) containing \(A\) such that \(f(a) \Linind_A B\), we have
\(\tp(B/Af(a)) \vdash \tp(B/Aa)\).

\item $p$ is residually dominated if it is so via any any pro-\(A\)-definable function $f$ sending $x$ to \(\Lin_A(x)\).

\item a type $q$ over $M$ is residually dominated over \(A\) if it is \(A\)-invariant
and $q|A$ is residually dominated.
\end{itemize}
\end{definition}

\begin{remark}
\label{models suffice resdom}
Similarly to domination (\cref{models suffice dom}), to check that a type $p$ over $A$ is residually dominated via \(f\) it suffices to show that for some realization $a\models p$ and some sufficiently saturated and homogeneous \(N\subsel M\), if \(f(a) \Linind_{A} N\) implies
\(\tp(N/Af(a)) \vdash \tp(N/Aa)\).
\end{remark}

\begin{remark}
\label{resdom alt}
\begin{enumerate}
\item If \(\k = \alg{\k}\) is a pure field, then a type over \(A\) is residually dominated if and only if it is dominated by \(\Lin_A\). Indeed, as \(\Lin_A\) is \(\Lin_A(\acleq(A))\)-definable, \cref{fork st emb} applied above \(\Lin_A(\acleq(A))\) provides a translation between forking in \(M\) and forking in \(\Lin_A = \alg{\k}\tensor\Lin_A\).

\item If \(A\subsel M\), then \(\Lin_A(M)\) is interdefinable over \(A\) with
\(\dcleq(\k(M))\). It follows that \(\tp(a/A)\) is residually dominated if and
only if for every \(B\subseteq\eq{M}\) containing \(A\) with \(\k(Aa)\)
algebraically independent from \(\k(B)\) over \(\k(A)\), we have
\(\tp(B/A\k(Aa))\vdash \tp(B/Aa)\); thus recovering earlier definitions
(\cite{has-ealy-mar,has-ealy-simon,Vic-ResDom}).
\end{enumerate}
\end{remark}

\begin{lemma}
\label{resdom going up}
Let \(B\subseteq \eq{M}\) contain \(A\). If \(\tp(a/A)\) is residually dominated
via some pro-\(A\)-definable \(f\) and if \(f(a) \Linind_{A} B\),
then \(\tp(a/B)\) is residually dominated via \(f\).
\end{lemma}

\begin{proof}
Let \(C\supseteq B\) and let us assume that \(f(a)\Linind_{B} C\). Since
\(\alg{\k}\tensor\Lin_A\) is stable, by transitivity, it follows that
\(f(a)\Linind_{A} C\).  By residual domination, we have
\(\tp(C/Af(a))\vdash \tp(C/Aa)\) and hence \(\tp(C/Bf(a))\vdash \tp(C/Ba)\).
\end{proof}



\begin{lemma}
\label{resdom orth}
If \(\tp(a/A)\) is residually dominated then \(\eq{\Gamma}(Aa) =
\eq{\Gamma}(A)\).
\end{lemma}

\begin{proof}
Let \(\gamma\in \eq{\Gamma}(Aa)\). Since \(\k\) and \(\Gamma\) are orthogonal,
we have \(\Lin_A(\acleq(A)\gamma) = \Lin_A(\acleq(A))\). It follows that \(\Lin_A(a)\Linind_A\gamma\). By residual domination, we have
\begin{equation*}
\tp(\gamma/A
\Lin_A(a))\vdash \tp(\gamma/Aa)
\end{equation*}
which has a unique realization. It follows
that \(\gamma\in\eq{\Gamma}(A \Lin_A(a)) = \Gamma(A)\), again by orthogonality.
\end{proof}

\begin{lemma}
\label{resdom push}
Assume that \(\Lin_A(\acleq(A))\subseteq A\). If \(\tp(a/A)\) is residually
dominated and \(c\in\dcleq(Aa)\), then \(\tp(c/A)\) is residually dominated.
\end{lemma}

We follow the proof of \cite[Proposition~3.32]{HHM}

\begin{proof}
Let \(B\subseteq M\) be such that \(\Lin_A(c)\Linind_A B\). Let
\(B'\equiv_{A\Lin_A(c)} B\), we have to show that \(B'\equiv_{Ac} B\). Let
\(d\equiv_{Ac} a\) be such that \(\Lin_A(d)\Linind_{A\Lin_A(c)} BB'\). By
transitivity, we have \(\Lin_A(d)\Linind_{A} B\) and \(\Lin_A(d)\Linind_{A} B'\)
and thus, by stationarity (see \cref{char Linind}.1), we have
\(B\equiv_{A\Lin_A(d)} B'\). By residual domination, it follows that
\(B\equiv_{Ad} B'\) and therefore \(B\equiv_{Ac} B'\).
\end{proof}

\section{Algebraically closed residue fields and dense value groups}

In this section we assume that \(M\) is a \(\Gamma\)-expansion of a henselian
valued field with algebraically closed characteristic zero residue field and
dense value group satisfying property \D. Let \(A\subseteq \eq{M}\).

\begin{lemma} Let \(a\) be a tuple in \(\K(\MM)\) such that \(\val(M(a)) =
\val(M)\).
\label{max model orth}
\begin{enumerate}
\item We have \(\tp_0(a/M)\vdash \tp(a/M)\).
\item If \(M\) is maximally complete, then \(\tp(a/M)\) is dominated by \(\k\).
\end{enumerate}
\end{lemma}

\begin{proof}
As \(\val(M(a)) = \val(M)\), the first statement follows from quantifier
elimination (\cref{EQ unramified}). Let us now prove that \(\tp(a/M)\) is
dominated by \(\k\). Let \(M\subsel N\subsel \MM\) be such that \(\k(Ma) \ind_M
N\). Then, as \(\k(M)\leq \k(Ma)\) is regular, \(\k(Ma)\) is linearly disjoint
from \(\k(N)\) over \(\k(M)\). As \(\Gamma(Ma) = \Gamma(M)\), by
\cite[Corollary~12.12]{HHM}, we have
\[\tp(N/M\k(Ma)) \vdash \tp_0(N/M\k(Ma)) \vdash \tp_0(N/Ma) \vdash \tp(N/Ma),\]
where the last implication follows from the first statement and \cref{type
game}.
\end{proof}

\begin{theorem}
\label{equiv epsilon}
Let \(a\) be a tuple from \(\K(\MM)\) such that \(p = \tp(a/M)\) is
\(A\)-invariant. The following statements are equivalent:
\begin{enumerate}[(i)]
\item \(p\) is dominated by \(\kInt\) over \(A\);
\item \(p\) is residually dominated over \(A\);
\item \(p\) is orthogonal to \(\Gamma\).
\end{enumerate}
\end{theorem}

\begin{proof}\phantom{}
\begin{description}[leftmargin=0pt,labelindent=\parindent,font=\rmfamily]
\item[(i)\imp(ii)] 
Assume \(\tp(a/A)\) is dominated by \(\kInt_A\). Recall that, as \(\k(M) =
\alg{\k(M)}\), residual domination over \(A\) is equivalent to domination by
\(\Lin_A\) (see \cref{resdom alt}). By \cref{dom fork}, \(\restr{p}{A}\) is
stationary, as \(p\) is \(A\)-invariant, and it suffices
to check that \(\Lin_A(a)\ind_A M\) implies \(a\ind_A M\).

Now, assume \(\Lin_A(a)\ind_A M\). By \cref{char kint}, we have \(\kInt_A(c) \subseteq
\acleq(A\Lin_A(a))\). By \cref{forking acl}, it follows that \(\kInt_A(a)
\subseteq \acleq(\Lin_A(a))\ind_A M\). As \(\tp(a/A)\) is dominated by
\(\kInt_A\), by \cref{dom fork}.2, we have \(a\ind_A M\).

\item[(ii)\imp(iii)] If \(\tp(a/A)\) is residually dominated, then \(p\) is
residually dominated (\emph{c.f.} \cref{resdom going up}). By \cref{resdom orth}, it
follows that \(p\) is orthogonal to \(\Vg\).

\item[(iii)\imp(i)] Assume that \(p\) is orthogonal to \(\Gamma\). We may assume
that \(M\) is maximally complete. By \cref{max model orth}, the type
\(\tp(a/M)\) is dominated by \(\k\). By \cref{descent intdom}, the type
\(\tp(a/A)\) is dominated by \(\kInt_A\).\qedhere
\end{description}
\end{proof}

\begin{remark}
The above statements are also equivalent to \(p\) being generically stable in this setting, but
we will further discuss this in \cref{gen st} in a wider context.
\end{remark}

\section{Between \texorpdfstring{\(M\)}{M} and \texorpdfstring{\(\ur{M}\)}{Mur}}

Through this section we let \(A \subseteq \eq{M}\) contain \(\G(\acleq(A))\), let \(a\) be a tuple in
\(\K(\MM)\) and let \(\epsilon \in \{1,1'\}\) (see \cref{geom sorts} for the definition of $\mathcal{G}$). We assume that $p_{0}=\tp_{0}(a/ M)$ is $\L(A)$-definable. 
\begin{definition}
Let \(\rho_{A,1}(a)\) enumerate \(\Lin_A(\MM) \cap \dcleq_1(\G(A)a)\).

In this section we will work with three different structures $M$, $M^{ur}$ and $M^{a}$ so we want to remind the reader about the notation introduced in \cref{setup}.
\end{definition}

\begin{lemma}
\label{compare Lin}
We have
\(\Lin_{\G(A),\epsilon}(a)\subseteq\dcleq_\epsilon(\rho_{A,1}(a))\).
\end{lemma}

\begin{proof}
Let us first assume that \(\epsilon = 1\). By \cite[Lemma~5.8]{RKVic-EIAKE} we have
\(\G_1(\dcleq_1(\MM))\subseteq \G(\MM)\). Let \(R\in \Gr_1(\acleq_1(\G(A)))\) and
\(e\in\red_1(R)(\dcleq_1(\G(A)a))\). Then
\[R\in \Gr_1(\acleq_1(\G(A)))\cap\dcleq_1(\MM)\subseteq \Gr(\acleq(A))\subseteq
\Gr(A).\] Similarly, \(e\in \G_1(\dcleq_1(\MM))\subseteq \MM\) and hence
\(e\in\red(R)(\MM)\cap\dcleq_1(\G(A)a) \subseteq\rho_{A,1}(a)\).

We now assume \(\epsilon = 1'\). Let \(R \in \Gr_{1'}(\acleq_{1'}(\G(A)))\) and
\(e\in\red_{1'}(R)(\dcleq_{1'}(\G(A)a))\). Note that \(R\) and \(e\) are in
\(\dcleq_{1'}(\MM)\). By \cite[Lemma~5.11]{RKVic-EIAKE}, there exits \(Q
\in\Gr_1(\MM_1)\) and \(\eta \in \red_1(Q)(\MM_1)\) such that:
\begin{itemize}
\item \(R\in \dcleq_{1'}(Q)\) and \(Q\) is definable from \(R\) in the pair
\((M_1,M_{1'})\);
\item \(e\in \dcleq_{1'}(\eta)\) and \(\eta\) is definable from \(e\) in the pair
\((M_1,M_{1'})\).
\end{itemize}
Given \(\sigma\in\Aut(\MM_1/\G(A))\), we can extend it to \(\MM_{1'}\). As \(R\)
has a finite orbit over \(\G(A)\) in \(\MM_{1'}\), it follows that
\(Q\in\acleq_1(\G(A))\). Similarly,
\(\eta\in\red_1(Q)(\dcleq_1(\G(A)a))\subseteq \Lin_{\G(A),1}(a)\subseteq
\rho_{A,1}(a)\). So \(e\in\dcleq_{1'}(\eta)\subseteq
\dcleq_{1'}(\rho_{A,1}(a))\). 
\end{proof}

\begin{remark}
    Note that the previous statement does not require $p_0$ to be $\L(A)$-definable. 
\end{remark}

\begin{lemma}
\label{code germ open ball} 
Let $b(a)$ be an $\L_{0}(M a)$-definable open ball with radius in
$\Vg(\dcleq_{1'}(\G(A)a))$. Then \(\germ{b}{p_0} \in \dcleq_\epsilon(\G(A),\Lin_A(M))\).
\end{lemma}

This is essentially \cite[Lemmas~5.9 and 5.13]{RKVic-EIAKE}, but we include a
proof for the sake of completeness.

\begin{proof}
By \cite[Proposition~4.3]{RKVic-EIAKE}, $\tp_{0}(a/ M_{1})$ is 
the
unique extension of $\tp_{0}(a/ M)$ and it is
$\L_{1}(\cG(A))$-definable. Given \cref{compare Lin}, we may assume
that $M=M_{1}$ --- growing it we may still assume that it suﬃciently saturated
and homogeneous. By \cite[Proposition~5.13]{RKVic-EIAKE}, the type $\tp_{0}(a/
M_{1'})$ is $\L_{1'}(\cG(A))$-definable. By
\cite[Lemma~5.9]{RKVic-EIAKE} applied in $M_{1'}$ over \(\G(A)\), the germ
$\germ{b}{p_{0}}$ is coded in $\cG_{1'}(\acleq_{1'}(\cG(A)))\cup
\Lin_{\cG(A),1'}(M_{1'})$ over $\cG(A)$ (in the structure $M_{1'}$). As
\(\germ{b}{p_{0}} \in \dcleq_{1'}(M_1)\) and 
\[
\K(\acl_{1'}(\G(A)))\cap\dcleq_{1'}(M)
= \K(\acleq(A)) \subseteq \G(A),
\]
by \cref{compare Lin}, we have \(\germ{b}{p_0}
\in \dcleq_{1'}(\G(A),\Lin_A(M))\). As any \(\sigma\in\Aut(M_1)\) extends to
\(M_{1'}\), we also have \(\germ{b}{p_0} \in \dcleq_{1}(\G(A),\Lin_A(M))\).
\end{proof}

\begin{proposition}
\label{resdom control}
Assume that \(\val(M(a))=\val(M)\). Then
\[\rv(M(a))\subseteq\dcleq_{\epsilon}(\cG(A),\rv(M),\Lin_A(M),\rho_{A,1}(a)).\]
\end{proposition}

By induction on an enumeration of \(a\), this is a consequence of the following
relative dimension one case:

\begin{lemma}
\label{control rel}
Let \(c\in\K(\MM)\) be such that \(q_0 = \tp_0(ac/M)\) is \(\L(A)\)-definable and that \(\val(M(ac)) = \val(M)\). Then \(\rv(M(ac)) \subseteq
\dcleq_{\epsilon}(\rv(M(a)),\cG(A),\Lin_A(M), \rho_{A,1}(ac))\).
\end{lemma}

\begin{proof}
By \cref{orth Vg ext ACVF}, we have $\tp_0(ac/M) \vdash\tp_{0}(ac/M_{0})$ and $\tp_{0}(ac/M_{0})$ is $\L_{0}(\G_0(A))$-definable and orthogonal to \(\Vg\) ---  in particular, it is stably dominated and germs of \(\L_0(M_0)\)-definable functions over \(q_0\) are strong (\emph{c.f.} \cref{strong germs}). We proceed by cases.\medskip

\begin{description}[leftmargin=0pt,labelindent=\parindent,font=\rmfamily]
\item[Case~1:] There is an \(\L_0(Ma)\)-definable finite union of closed balls \(B(a)\) containing \(c\) which is minimal among \(\L_0(Ma)\)-definable finite union of balls.

\begin{claim}
We have \(B(a) \in \dcleq_0(\G_0(A)a)\).
\end{claim}

\begin{claimproof}
Note first that \(c\) is generic in \(B(a)\) over \(M_0a\). Indeed, if \(B'(a)\)
is an \(\L_0(M_0a)\)-definable finite union of closed balls containing
\(c\), then the union of its \(Ma\)-conjugates must contain \(B(a)\). But, since
$\tp_0(ac/M) \vdash\tp_{0}(ac/M_{0})$, any of these conjugates must contain
\(B(a)\).

Then \(\germ{B}{p_{0}}\) is coded (in \(M_0\)) by \(\ulcorner q_0 \urcorner\) over
\(\cG_0(A)\). So \(\germ{B}{p_0}\in\dcleq_0(\cG_{0}(A))\). As germs over \(p_0\) are
strong, it follows that \(B(a)\in\dcleq_0(\cG_{0}(A)a)\). \end{claimproof}\medskip

Let $b\in B(a)$ be the ball containing $c$. By \cite[Lemma~4.3.4]{HilRK-EIAKE}
there is an $\L_{0}(\code{B(a)})$-definable injective map $B(a) \to \RV^{n}$. So
there is a tuple $\nu(ac)$ in $\RV(\dcleq_{0}(\G_0(A)ac))$ such that $b$ is
$\L_0(\G_{0}(A)\nu(ac)a)$-definable. By \cite[Lemma~2.6]{RKVic-EIAKE}, there is
an $\L_{0}(M a)$-definable finite set $H(a)\subseteq \K$ such that $H(a) \cap
b=\{h(ac)\}$ is a singleton. By \cite[Lemma~4.3.13]{HilRK-EIAKE} we have:
\[\rv(M(ac)) \subseteq \dcleq_0(\rv(M(a)), \nu(ac), \rv(c-h(ac))).\]

Let \(\eta(ac) = \rv(c-h(ac))\). To conclude case 1, there remains to show the
following:

\begin{claim}
\label{germ rv gen}
We have \(\nu(ac),\eta(ac)\in\dcleq_\epsilon(\G(A),\Lin_A(M),\rho_{A,1}(ac))\).
\end{claim}

\begin{claimproof}
By minimality of \(B(a)\), note that \(\val(c-h(ac))\) is equal to the radius of
\(b\) which thus lies in \(\dcleq_0(\G_0(A))\cap\Vg(\MM)\), as \(\tp_0(a/M_0)\) is
orthogonal to \(\Vg\). By \cref{code germ open ball}, it follows that
\(\germ{\eta}{q_0}\in \dcleq_{\epsilon}(\cG(A),\Lin_{A}(M))\). As germs over
\(q_0\) are strong, we have \(\eta(ac) \in
\dcleq_{0}(\cG_{0}(A),\germ{\eta}{q_{0}}, ac) \subseteq
\dcleq_\epsilon(\cG(A),\Lin_A(M),ac)\). As \(\eta(ac)\) lies in
\(\alg{\k}\tensor\Lin_{\G(A)}\) which is stably embedded and eliminates
imaginaries (\cref{LinA EI}) it follows that
\[\eta(ac)\in\dcleq_\epsilon(\cG(A),\Lin_A(M),\Lin_{\G(A),\epsilon}(ac))\subseteq
\dcleq_\epsilon(\cG(A),\Lin_A(M),\rho_{A,1}(ac)),\] where the last inclusion
follows from \cref{compare Lin}.

The proof for \(\nu(ac)\) is similar as \(\val(\nu(ac)) \in \dcleq_0(\G_0(A))\cap\Vg(\MM)\) by orthogonality, and hence \(\nu(ac) \in \Lin_{\G(A),\epsilon}\).
\end{claimproof}\medskip

\item[Case~2:] No such \(B(a)\) exists.\medskip

By \cite[Lemma~4.3.10]{HilRK-EIAKE}, there is a pro-$\L_{0}(M)$-definable map
$\eta$ into a power of $\RV$ such that $\germ{\eta}{q_{0}} \in \dcleq(A)$ and
$\rv(M(ac)) \subseteq \dcleq_{0}(\rv(M(a)),\eta(ac))$. We now conclude as in
\cref{germ rv gen}.\qedhere
\end{description}
\end{proof}

\begin{corollary}
\label{resdom from epsilon}
If \(\tp_\epsilon(a/\G(A))\) is residually dominated in \(M_\epsilon\), then
\(\tp(a/A)\) is residually dominated via \(\rho_{A,1}\) in \(M\).
\end{corollary}

\begin{proof}
We may assume that \(\rho_{A,1}(a)\Linind_{\G(A)}M\). Since
\(\Lin_{\G(A),\epsilon}(a)\subseteq\dcleq_\epsilon(\rho_{A,1}(a))\), we have 
\begin{align*}
\tp(a/A)\cup\tp(\rho_{A,1}(a)/M)
&\vdash \tp_\epsilon(a/\G(A))\cup \tp_\epsilon(\Lin_{\G(A),\epsilon}(a)/M)\\
&\vdash \tp_\epsilon(a/M)\hspace{20pt} \text{(by residual domination)}\\
&\vdash \tp_0(a/M).
\end{align*}
Moreover, by \cref{resdom control}, we have
\(\tp(a/A)\cup\tp(\rho_{A,1}(a)/M)\vdash \tp(\rv(M(a)/M))\) and, by \cref{EQ
RV}, we have \(\tp_0(a/M)\cup\tp(\rv(M(a))/M) \vdash \tp(a/M)\). So
\begin{align*}
\tp(a/A)\cup\tp(\rho_{A,1}(a)/M)
&\vdash \tp_0(a/M)\cup\tp(\rv(M(a))/M)\\
&\vdash\tp(a/M). \tag*{\qedhere}
\end{align*}
\end{proof}

\section{The main equivalence}

Let \(A \subseteq \eq{M}\) contain \(\G(\acleq(A))\) and let \(a\) be a tuple in
\(\K(\MM)\).

\begin{theorem}
\label{equiv resdom}
Let \(\epsilon\in\{1,1'\}\), the following statements are equivalent:
\begin{enumerate}[(i)]
\item \(\tp(a/A)\) is residually dominated;

\item for some \(c\equiv_A a\) in \(\MM\) such that \(\tp_0(c/M)\) is
\(\L(\acleq(A))\)-invariant, we have \(\val(M(c)) = \val(M)\);

\item \(\tp_{\epsilon}(a/\G(A))\) is residually dominated.
\end{enumerate}
Furthermore, for \(\eta\in\{0,1,1'\}\), the type \(\tp_\eta(a/\G_0(A))\) is then
stationary and its non forking extension is
\(\L_\eta(\G_0(A))\)-definable, orthogonal to \(\Gamma\) and consistent with
\(\tp(a/A)\).
\end{theorem}

\begin{proof} We first prove the equivalence.
\begin{description}[leftmargin=0pt,labelindent=\parindent,font=\rmfamily]
\item[(i)\imp(ii)] 
By \cite[Theorem~4.1]{RKVic-EIAKE} and compactness, we find \(c\equiv_A a\) in
\(\MM\) such that the type \(\tp_1(c\Lin_A(c)/M)\) is
\(\L(\acleq(A))\)-invariant. By \cref{char Linind}.2, it follows that
\(\Lin_A(c)\Linind_{\G(A)} M\). Hence, by \cref{resdom going up}, the type
\(\tp(c/M)\) is residually dominated and, by \cref{resdom orth}, we have
\(\val(M(c)) = \val(M)\).

\item[(ii)\imp(iii)] Let \(c\) be as in (ii). As \(\G_0(\acleq(A))\subseteq A\),
by \cref{orth Vg ext ACVF}.2, the type \(\tp_{1'}(c/M_1')\) is
\(\L_{1'}(\G_0(A))\)-definable and orthogonal to \(\Gamma\). In particular, it
is \(\L_{1'}(\G(A))\)-definable. It follows from \cref{equiv epsilon} that
\(\tp_{1'}(a/\G(A)) = \tp_{1'}(c/\G(A))\) is residually dominated --- so is
\(\tp_{1}(a/\G(A))\) by \cref{resdom from epsilon} applied in
\(M_1\).

\item[(iii)\imp (i)] This is \cref{resdom from epsilon}.
\end{description}

Let us now assume that all three statement hold. Let \(c\) be as in (ii). By
\cref{orth Vg ext ACVF}, the type \(\tp_\eta(c/M_\eta)\) is
\(\L_\eta(\G_0(A))\)-definable and orthogonal to \(\Gamma\) and, as shown above,
it is residually dominated. As \(\Lin_{\G_0(A),\eta}\) is stable and
\(\tp_\eta(\Lin_{\G_0(A),\eta}(c)/M_\eta)\) is \(\L_\eta(\G_0(A))\)-definable,
it is stationary over \(\G_0(A)\) and so is \(\tp_\eta(c/M_\eta)\) by \cref{dom
fork}.1.
\end{proof}

\begin{remark}
\label{equiv orth}
Let us assume that \(\tp(a/M)\) is \(\L(A)\)-invariant.
\begin{itemize}
\item If \(\tp(a/A)\) is residually dominated, then by \cref{resdom going up},
\(\tp(a/M)\) is residually dominated and hence orthogonal to \(\Vg\) by
\cref{resdom orth}.
\item Conversely, if \(\tp(a/M)\) is orthogonal to \(\Vg\), then \(\val(M(a)) =
\val(M)\). As \(\tp_0(a/M)\) is \(\L(A)\)-invariant, it follows that
\(\tp(a/A)\) is residually dominated, by \cref{equiv resdom}.
\end{itemize}
\end{remark}

\begin{remark}
\label{strong stat}
If \(\tp(a/A)\) is residually dominated, then it has a unique
\(\L(A)\)-invariant extension to an \(\L_0(M)\)-type --- in other words, the
type of any \(c\) as in (ii) is unique. Indeed, if \(c\equiv_A a\) in \(\MM\) is
such that \(\tp_0(c/M)\) is \(\L(A)\)-invariant, then \(\tp_0(c/M_1)\) is
\(\L_1(\G(A))\)-definable; proceed as in \cref{orth Vg ext ACVF}. By Property \D
and compactness, we may assume that \(q_1 = \tp_1(c/M_1)\) is
\(\L_1(\G(A)\cup\eq{\Vg}(\acleq(A)))\)-definable (\cite[ Lemma~4.7]{RKVic-EIAKE}).
Recall that the unique non forking extension \(p_1\) of \(\tp_1(a/\G_0(A))\) is
consistent with \(\tp(a/A)\) and hence contains \(\tp_1(a/\G(A)\cup
\eq{\Vg}(\acleq(A)))\). It follows that \(q_1 = p_1\).
\end{remark}

\begin{remark}
The theorem also holds for \(a\in\G_0(\MM)\), replacing \(\val(M(a))\) by
\(\dcleq_0(Ma)\cap\Gamma(\MM)\) in statement (ii). The implication (i)\imp (ii) is
proved exactly as above. The implication (iii)\imp (ii) is proved in a similar
manner.

Now assume (ii) holds, \emph{e.i.} that \(\tp(a/M)\) is
\(\L(\acleq(A))\)-invariant and \(\dcleq_0(Ma)\cap\Gamma(M) = \val(M)\). There
is some \(\L_0(a)\)-definable type \(p_a\), orthogonal to \(\Gamma\),
concentrating on some \(\K^m\times\k^n\) and finitely satisfiable in \(\MM\),
such that for any \(c\models \restr{p_a}{Ma}\), we have \(a\in\dcleq_0(c)\) (see
the proof of \cite[Lemma~4.4.4]{HilRK-EIAKE}). Let \(c_{\K} = c\cap \K\). We have
\(\val(Mc_\K) \subseteq \dcleq_0(Mac_\K)\cap\Gamma(\MM) =
\dcleq_0(Ma)\cap\Gamma(\MM) = \val(M)\). It follows, by \cref{equiv resdom}, that
\(\tp(c_\K/A)\) is residually dominated, but then so is \(\tp(c/A)\). By
\cref{resdom push}, the type \(\tp(a/A)\) is also residually dominated.
Similarly, we see that \(\tp_\epsilon(a/A)\) is residually dominated.

The equivalence between (i) and (iii) for \(a\in\G(M)\) seems harder as it is
not obvious how to build good resolutions in that context.
\end{remark}

From \cref{equiv resdom}, we can deduce base change for residual domination along invariant extensions:

\begin{corollary}
\label{base change}
Assume that \(\tp_0(a/M)\) is \(\L(A)\)-invariant. For every \(B\subseteq
\eq{M}\) containing \(\G(\acleq(AB))\), the following statements are equivalent:
\begin{enumerate}[(i)]
\item \(\tp(a/B)\) is residually dominated;
\item \(\tp(a/A)\) is residually dominated.
\end{enumerate}
\end{corollary}

\begin{proof}
Let us assume that \(\tp(a/B)\) is residually dominated. By \cref{equiv resdom}
and \cref{strong stat}, the type \(\tp_0(a/M)\) is then orthogonal to \(\Gamma\)
and hence, again by \cref{equiv resdom}, the type \(\tp(a/A)\) is residually
dominated. The converse is symmetric.
\end{proof}

\section{Valued fields with operators}
\label{operator}

Let \(M_\nabla\) be a sufficiently saturated and homogenous expansion of \(M\)
is some language \(\L_\nabla\) with stably embedded orthogonal residue field and
value group --- indices \(_\nabla\) indicate that we consider the structure
\(M_\nabla\). Let \(\nabla\) be a pro-\(\L_\nabla\)-definable map from \(\K\) into a
power of \(\K\). We assume that
\begin{equation}\label{property H}\tag{$\mathrm{H}$}
\text{  for every tuple \(a\) in \(\K(M)\), \text{we have} \(\tp(\nabla(a)) \vdash
\tp_{\nabla}(a)\).} 
\end{equation}
We also assume that \(M_\nabla\) admits a spherically complete elementary
extension.

\begin{remark}
\label{operator examples}
This covers various theories of valued fields with operators considered in the literature.
\begin{enumerate}
\item If \(M_\nabla\) is \(\sigma\)-henselian equicharacteristic zero valued
difference field, then \eqref{property H} holds for \(\nabla(x) =
(\sigma^i(x))_{i\geq 0}\) --- see \cite{diffRV}. If we further assume that the
automorphism is multiplicative (\cite{Pal}, see also
\cite[Section~2.4]{HilRK-EIAKE}), then \(\k\) and \(\Gamma\) are stably embedded
and orthogonal. Moreover the induced structure \(\k\) is a reduct of a real closed field and hence is \(o\)-minimal.
\item If \(M_\nabla\) is an equicharacteristic zero \(D\)-henselian monotone
differential valued field, then \(M_\nabla\) verifies all required properties
with \(\nabla(x) = (\partial^i(x))_{i\geq 0}\) --- see \cite{Scanlon}.
\item This is also the case if \(M_\nabla\) is an equicharacteristic zero henselian
field with a generic derivation --- see \cite{cubides-point}.
\end{enumerate}
\end{remark}

Let \(\MM_\nabla\) be a sufficiently saturated and homogeneous elementary extension
of \(M_\nabla\). We may assume that \(\MM\) is the \(\L\)-reduct of \(\MM_\nabla\).

\begin{theorem}
\label{equiv resdom op}
Let \(A\subseteq\eq{M}_\nabla\) contain \(\G(\acleq_\nabla(A))\) and let \(a\) be a tuple
in \(\K(\MM)\). The following are equivalent:
\begin{enumerate}[(i)]
\item the type \(\tp_\nabla(a/A)\) is residually dominated and for some
\(c\equiv_{A,\nabla} a\) in \(\MM\), the type \(\tp_1(\nabla(c)\Lin_A(c)/M)\) is
\(\L_\nabla(\acleq_\nabla(A))\)-invariant;

\item for some \(c\equiv_{A,\nabla} a\) in \(\MM\), the type
\(\tp_0(\nabla(c)/M)\) is \(\L_\nabla(\acleq_\nabla(A))\)-invariant and we have
\(\val(M(\nabla(c))) = \val(M)\);

\item the type \(\tp(\nabla(a)/\G(A))\) is residually dominated in \(\MM\).
\end{enumerate}
\end{theorem}

\begin{proof}\phantom{}
\begin{description}[leftmargin=0pt,labelindent=\parindent,font=\rmfamily]
\item[(i)\imp(ii)]
Let \(c\) be as in (i). By \cref{char Linind}.2, we have
\(\Lin_{A,\nabla}(c) \Linind_{\G(A)} M\). Thus, by \cref{resdom going up}, the
type \(\tp(\nabla(c)/M)\) is residually dominated in $\MM$ and, by \cref{resdom orth}, we have
\(\val(M(\nabla(c))) = \val(M)\). 

\item[(ii)\imp(iii)] Let \(c\) be as in (ii). Since \(\tp_0(\nabla(c)/M)\) is
\(\L_\nabla(\acleq_\nabla(A))\)-invariant, by \cref{orth Vg ext ACVF} applied in
\(M_\nabla\), it is \(\L(\G_0(A))\)-definable. It follows from \cref{equiv
resdom} that \(\tp(\nabla(c)/\G(A)) = \tp(\nabla(a)/\G(A))\) is residually
dominated in $\MM$.

\item[(iii)\imp(i)] Let us assume that \(\tp(\nabla(a)/\G(A))\) is residually
dominated in $\MM$. By \cref{equiv resdom}, the type \(\tp_1(\nabla(a)/\G(A))\) is residually
dominated in $\MM_{1}$.

\begin{claim}
The type \(\tp_1(\nabla(a)\Lin_A(a)/\G(A))\) is residually dominated in $\MM_{1}$.
\end{claim}

\begin{claimproof}
Let \(c\) enumerate \(\Lin_A(a)\). By elimination of imaginaries (\emph{c.f.} \cref{LinA
EI}), \(\Lin_{\G(A),1}(\nabla(a)c)\) in inter-\(\L_1\)-definable with
\(c\Lin_{\G(A),1}(\nabla(a))\).  If \(B\subseteq \eq{\MM_1}\) contains \(\G(A)\) and
is such that \(\Lin_{\G(A),1}(\nabla(a)c)\Linind_{\G(A)} B\), then
\(\Lin_{\G(A),1}(\nabla(a))\Linind_{\G(A)} Bc\). By residual domination, we have
\(\tp_1(Bc/\G(A)\Lin_{\G(A),1}(\nabla(a)))\vdash \tp_1(Bc/\G(A)\nabla(a))\) and hence
\begin{align*}\
\tp_1(B/\Lin_{\G(A),1}(\nabla(a)c))
&\vdash \tp_1(B/\G(A)\Lin_{\G(A),1}(\nabla(a)),c)\\
&\vdash \tp_1(B/\G(A)\nabla(a)c).\tag*{\claimqed}
\end{align*}
\end{claimproof}

Then \(\tp_1(\nabla(a)\Lin_A(a)/\G(A))\) is stationary and its non forking
extension is consistent with \(\tp(\nabla(a)\Lin_A(a)/\G(A))\) and
\(\L(\G(A))\)-definable --- see \cref{char Linind}.1 and \cref{dom fork}.1.

There remains to show that \(\tp_{\nabla}(a/A)\) is residually dominated in $\MM_{\nabla}$. Let us assume
that \(\Lin_{A,\nabla}(a)\Linind_A M\). As \(\Lin_{\G(A)}(\nabla(a)) \subseteq
\Lin_{A,\nabla}(a)\), we have \(\Lin_A(\nabla(a))\Linind_{\G(A)} M\). As
\(\tp(\nabla(a)/\G(A))\) is residually dominated in $\MM$, we have
\begin{align*}
\tp_\nabla(a/\G(A))\cup\tp(\Lin_{A,\nabla}(a)/M)
&\vdash\tp(\nabla(a)/\G(A))\cup \tp(\Lin_{\G(A)}(\nabla(a))/M)\\
&\vdash\tp(\nabla(a)/M)\\
&\vdash\tp_{\nabla}(a/M).\tag*{\qedhere}
\end{align*}
\end{description}
\end{proof}

\begin{remark}
As in \cref{equiv orth}, if \(\tp(a/M)\) is \(\L_\nabla(A)\)-invariant, then
\(\tp(a/A)\) is residually dominated if and only if \(\tp(a/M)\) is orthogonal
to \(\Gamma\). Moreover, note that by \cref{equiv resdom}, the type \(\tp(\nabla(a)/\G(A))\) is residually dominated if and only if the type \(\tp_1(\nabla(a)/\G(A))\) is residually dominated. 
\end{remark}

\section{Lifting tameness}

In this section we show that some model-theoretic tameness of the residue field
is inherited by types that are orthogonal to the value group. The two main
results of this section are \cref{gen st equiv}, for generic stability, and
\cref{gen simp equiv}, for generic simplicity. 

Let \(M_\nabla\) be as in section \cref{operator}. Note that the case \(M_\nabla
= M\) is covered by taking \(\nabla\) to be the identity function.

\subsection{Generic stability}
\label{gen st}

Recall the definition of generic stability (see \cref{def:gen st}).

\begin{proposition}
\label{gen st equiv}
Assume that \(\k_\nabla\) is stable. Let \(A \subseteq\eq{M}_\nabla\) contain
\(\G(\acleq_\nabla(A))\cup\eq{\Lin_{A,\nabla}}(\acleq_\nabla(A))\) and let \(a\) be a tuple in
\(\K(\MM)\). The following are equivalent:
\begin{enumerate}[(i)]
\item \(\tp_\nabla(a/A)\) is residually dominated;
\item \(\tp_\nabla(a/A)\) is stationary and its global non forking extension is
generically stable.
\end{enumerate}
\end{proposition}

\begin{proof}\phantom{}
\begin{description}[leftmargin=0pt,labelindent=\parindent,font=\rmfamily]
\item[(i)\imp(ii)]
Assume that \(\tp_\nabla(a/A)\) is residually dominated. As
\(\tp_\nabla(\Lin_A(a)/A)\) is stationary (see \cref{char Linind}), so is
\(\tp_\nabla(a/A)\). Its global non forking extension is generically stable by
\cref{resdom gen st}.
\item[(ii)\imp(i)] Assume that the unique non forking extension \(p\) of
\(\tp_\nabla(a/A)\) is generically stable. We may assume that \(a\models p\).
For any \(\gamma\in\val(M(\nabla(a)))\), the type \(\tp_\nabla(\gamma/M)\) is
generically stable in a totally ordered set. It follows that
\(\tp_\nabla(\gamma/M)\) is realized and hence \(\gamma\in \val(M)\). As
\(\tp_\nabla(a/M)\) is generically stable over \(A\), it is
\(\L_\nabla(A)\)-definable (see \cite[Proposition~2.1]{PilTan} or \cite[Proposition~9.6]{Casanovas}). By \cref{equiv resdom op}, it follows
that \(\tp_\nabla(a/A)\) is residually dominated.\qedhere
\end{description}
\end{proof}

\begin{remark}
\cref{gen st equiv} does not necessarily hold for imaginary elements. If
\(\Gamma\) is a pure ordered abelian group, then \(\Gamma/n\Gamma\) is stable.
If it is infinite, it contains non algebraic stable (and hence generically
stable) types. This also shows that, when \(\k\) is stable, the stable part over
\(A\) can be strictly larger than \(\acl(A)\cup\Lin_A\).
\end{remark}

\subsection{Generic simplicity}

Let \(M\) be some structure and \(A\subseteq\eq{M}\).

\begin{definition}
\label{def:gen simple}
A type \(p\) over \(A\) is \emph{generically simple} if, for every \(a\models
p\) and every tuple \(b\in\eq{M}\) with \(b\ind_A a\), we have \(a\ind_A b\).
\end{definition}

\begin{fact}[{\cite[Corollary~3.11 and Lemma~3.8]{Sim-GenSimp}}]
\label{prop gen simp}
Let \(a\) be a tuple in \(M\) be such that \(\tp(a/A)\) is generically simple.
\begin{enumerate}
\item Assume \(A\) is an extension base. If   \(B\subseteq \eq{M}\) contains \(A\)
is such that \(a\ind_A B\), then \(\tp(a/B)\) is generically simple.
\item If \(c\in\dcleq(Aa)\), then \(\tp(c/A)\) is generically simple.
\end{enumerate}
\end{fact}

Let \(M_\nabla\) (and \(M\)) be, once again, as in \cref{operator}.

\begin{proposition}
\label{gen simp equiv}
Assume that \(\k_\nabla\) is simple and that \(M_\nabla\) is \(\NTP_2\). Let
\(A\subseteq\eq{M}_\nabla\) be an extension base containing
\(\G(\acleq_\nabla(A))\). Let \(a\)  be a tuple in \(\K(\MM)\). The following
are equivalent:
\begin{enumerate}[(i)]
\item \(\tp_\nabla(a/A)\) is residually dominated;
\item \(\tp_\nabla(a/A)\) is generically simple and there exists
\(c\equiv_{A,\nabla} a\) such that \(c\ind^\nabla_A M\) and
\(\tp_0(\nabla(c)/M)\) is
\(\L_\nabla(\acleq_{\nabla}(A))\)-invariant.
\end{enumerate}
\end{proposition}

\begin{proof}\phantom{}
\begin{description}[leftmargin=0pt,labelindent=\parindent,font=\rmfamily]
\item[(i)\imp(ii)] Let us assume that \(\tp_\nabla(a/A)\) is residually
dominated. Let \(b\in\eq{M}\) be such that \(b\ind^\nabla_A a\). By \cref{fork st emb},
we have \(\eq{\Lin}_{A,\nabla}(b)\ind^\nabla_{\eq{\Lin}_{A,\nabla}(A)}
\eq{\Lin}_{A,\nabla}(a)\). As \(\k\) is simple, so is \(\eq{\Lin}_{A,\nabla}\)
and hence
\begin{equation}\label{a nf b}
\eq{\Lin}_{A,\nabla}(a)\ind^\nabla_{A} b.
\end{equation}
As \(\alg{\k}\tensor\Lin_A\) is stable and eliminates imaginaries, the type \(\tp_\nabla(\Lin_{A,\nabla}(a)/b)\) is consistent with an \(\L_1(\G(A))\)-definable type over \(M_1\). In particular, we have \(\Lin_{A,\nabla}(a)\Linind_{A}b\). By
residual domination, we have
\[\tp_\nabla(a/A)\cup\tp_\nabla(\Lin_{A,\nabla}(a)/Ab)\vdash \tp_\nabla(a/Ab).\]

If \(\tp_\nabla(a/Ab)\) divides over \(A\), there is a \(X\in\tp_\nabla(a/A)\),
an \(\L_\nabla(A)\)-definable map \(f\) into \(\Lin_{A,\nabla}\) and a \(Y_b\in
\tp_\nabla(f(a)/Ab)\) such that \(X \cap f^{-1}(Y_b)\) divides over \(A\) and
hence so does \(f(X)\cap Y_b\); but this contradicts \cref{a nf b}. So
\(\tp_\nabla(a/Ab)\) does not divide over \(A\) and neither does it fork over
\(A\) as \(M_\nabla\) is \(\NTP_2\) and \(A\) is an extension base (see
\cite[Theorem~1.2]{CheKap}). So  \(\tp_\nabla(a/A)\) is generically simple.

Moreover, Let \(c\equiv_{A,\nabla} a\) be such that \(c\ind_A M \). In
particular, we have \(\Lin_{A}(\nabla(c))\Linind_A M\) and thus
\(\tp_1(\nabla(c)/M)\) is the \(\L(\G_0(A))\)-definable extension of
\(\tp_1(\nabla(c)/\G(A))\) --- see \cref{equiv resdom}.

\item[(ii)\imp(i)] Let us assume that \(\tp_\nabla(a/A)\) is generically simple
and let \(c\) be as in (ii). Let \(\gamma\in\Vg(\acleq_{\nabla}(Ma))\). By
\cref{prop gen simp}, the type \(\tp_\nabla(\gamma/M)\) is generically simple.
But \(\Vg\) is ordered, so this type is realized. In particular, we have
\(\val(M(\nabla(a))) = \val(M)\). By \cref{equiv resdom op}, the type
\(\tp_\nabla(c/M)\) is therefore residually dominated.\qedhere
\end{description}
\end{proof}

\begin{remark} In many examples, for any tuple \(a\) in \(\K(\MM)\), there
exists \(c\equiv_{\acl(A),\nabla} a\) such that \(\tp_0(\nabla(c)/M)\) is
\(\L_\nabla(\acleq_\nabla(A))\)-invariant (see
\cite[Corollary~3.5.6]{HilRK-EIAKE} and \cite[Corollary~4.8]{RKVic-EIAKE}); and
in that case we have \(c\ind^\nabla_A M\) (\emph{e.g.} \cref{qf fork}). In other
words, the second part of condition (ii) always holds in these examples.
\end{remark}

\subsection{Extension bases}

Let us consider two algebraic examples of interest where \cref{gen simp equiv}
applies.

\begin{proposition}
\label{pseudolocal and VFA}
Consider the theory of either:
\begin{enumerate}[a)]
\item an equicharacteristic zero ultraproduct of \(p\)-adic fields;
\item an equicharacteristic zero existentially closed multiplicative difference
valued fields.
\end{enumerate}
Then the theory is \(\NTP_2\), the valued group satisfied property $\D$, the residue field is
simple and any set is an extension base.
\end{proposition}

Note that the second example includes the limit theory, as \(p\) goes to
infinity, of the Frobenius automorphism acting on an algebraically closed valued
field of characteristic $\p$, first considered in \cite{Hru-ACFA}.

\begin{proof}
Let \(M\) be an  equicharacteristic zero ultraproduct of \(p\)-adic fields. The
residue field is pseudo-finite and hence simple. 
 By \cite[Corollary $5$]{inpminimal} any non-principal ultraproduct of the $\p$-adics is inp-minimal, therefore $\NTP_{2}$. The value group is elementarily
equivalent to \(\ZZ\) so it satisfies Property $\D$. Finally any \(A\subseteq
\eq{M}\) is an extension base by \cite[Theorem~5.11]{Hos-ExtBase}.

Now, let \(M_\nabla\) be a sufficiently saturated and homogeneous existentially
closed multiplicative difference valued field. By
\cref{operator examples}.1, the value group is \(o\)-minimal and
therefore has Property \D. The residue field is a model of \(\ACFA\), the theory
of existentially closed difference fields, which is an unstable supersimple
theory and it eliminates imaginaries by results of Z. Chatzidakis and E. Hrushovski in \cite{ACFA}.
By
\cite[Theorem $4.1$]{VFANTP} the structure \(M_{\nabla}\) is $\NTP_{2}$.

There remains to show that any \(A = \acleq(A)\subseteq\eq{M}\) is an extension
base. We follow the strategy of \cite[Corollary~8.5]{DorHal}\footnote{We thank
Y. Dor and Y. Halevi for making us aware of this result and its proof.}. Let
\(a\) be a tuple in \(\K(\MM)\). By \cite[Theorem~3.1.1]{HilRK-EIAKE}, there
exists \(c\equiv_{\acleq_\nabla(A)} a\) such that \(\tp_0(\nabla(c)/M)\) is
\(\L_\nabla(\G_0(A))\)-definable. It now suffices to prove the following:

\begin{claim}
\label{qf fork}
The type \(\tp_\nabla(c/M)\) does not divide over \(A\).
\end{claim}

\begin{claimproof}
We first assume that \(A = N_\nabla\subsel M_\nabla\) is spherically complete.
Since \(\tp_0(\nabla(c)/M)\) does not fork over \(\G_0(A)\) and \(\k(A)\) is
algebraically closed, the residue field \(\k(\A(\nabla(c)))\) is linearly
disjoint from \(\k(M)\) over \(\k(\A)\). Also, we have
\(\Gamma(\A(\nabla(c)))\cap\Gamma(M) = \Gamma(\A)\). It follows that
\begin{align*}
\tp_0(\nabla(c)/\A) \cup \tp_{\nabla}(\k(\A(\nabla(c)))/\k(M))\cup
\tp_{\nabla}(\Gamma(\A(\nabla(c)))/\Gamma(M))
&\vdash
\tp(\nabla(c)/M)\\
&\vdash \tp_{\nabla}(c/M).
\end{align*}
Indeed, let \(d\) realize the left had side. Then, by orthogonality of \(\k\)
and \(\Gamma\), we may assume that \(\k(\A(\nabla(c))) = \k(\A(\nabla(d)))\) and
\(\Gamma(\A(\nabla(c))) = \Gamma(\A(\nabla(d)))\). By
\cite[Proposition~12.11]{HHM} and
its proof, the fields \(M(\nabla(c))\) and \(M(\nabla(d))\) are isomorphic over
\(\RV(M(\nabla(c))) = \RV(M(\nabla(d)))\). By field quantifier elimination
(\cref{EQ RV}), we have \(\tp(\nabla(c)/M) = \tp(\nabla(d)/M)\) and hence
\(\tp_\nabla(c/M) = \tp_\nabla(d/M)\).

As a definable set divides in \(\ACFA\) if and only its Zariski closure does,
the type \(\tp_{\nabla}(\k(\A(\nabla(c)))/\k(M))\) does not divides over
\(\k(\A)\). Also, as \(\Gamma\) eliminates quantifiers (as a difference ordered
group), the type \(\tp_{\nabla}(\Gamma(\A(\nabla(c)))/\Gamma(M))\) does not
divides over \(\Gamma(\A)\). It follows that \(\tp_\nabla(c/M)\) does not divide
over \(A\) (see \cite[Theorem~4.4]{Vic-ResDom}).

Now, if \(\A\) is any subsection of \(M\), let \(N_\nabla\subsel M_\nabla\)
contain \(A\). If \(\tp_\nabla(c/M)\) divides over \(A\), we find
\(X\in\tp_\nabla(c/M)\) and an \(\L_\nabla(A)\)-indiscernible sequence \(I\)
containing \(\code{X}\). Growing the sequence, we can find an
\(\L_\nabla(N)\)-indiscernible subsequence \(J\subseteq I\). Conjugating by an
automorphism, we may assume that \(\code{X}\in J\). So \(X\) also divides over
\(N\) and \(\tp_\nabla(a/M)\) divides over \(N\). This contradicts the first
case we considered.
\setcounter{claimqed}{0}
\end{claimproof}
\end{proof}

\sloppy
\printbibliography

\noindent\textsc{Pablo Cubides Kovacsics. Departamento de Matemáticas, Universidad de los Andes, Carrera 1 no. 18A - 12, Edificio H, Bogotá, Colombia.}

Email address:\textsf{p.cubideskovacsics@uniandes.edu.co}.\medskip

\noindent\textsc{Silvain Rideau-Kikuchi. DMA, École normale supérieure, Université PSL, CNRS, 75005 Paris, France.}

Email address: \textsf{silvain.rideau-kikuchi@ens.fr.}.\medskip

\noindent\textsc{Mariana Vicar\'ia. University of Chicago, Department of Mathematics, Eckhart Hall, 5734 S University Ave, Chicago, IL, 60637.}

Email address: \textsf{mvicaria@uchicago.edu}.
\end{document}